\documentstyle[a4,amscd,11pt]{amsart}

\newcommand{\mnote}[1]{}                        	
\newcommand{\closeup}{\setlength{\itemsep}{-2pt}}
\newcommand{\lars}{{\bf *LA*}}

\renewcommand{\Re}{{\Bbb R}}
\newcommand{\Ric}{\operatorname{Ric}}    	       
\newcommand{\ric}{\operatorname{ric}}           	
\newcommand{\tRic}{\widetilde{\Ric}}    	
\newcommand{\tomega}{\widetilde{\omega}}	    
\newcommand{\te}{\tilde{e}} 			   
\newcommand{\Za}{{\Bbb Z}}
\newcommand{\Co}{{\Bbb C}}
\newcommand{\n}{\nabla}			
\newcommand{\tn}{\widetilde{\nabla}}	
\newcommand{\hn}{\widehat{\nabla}}	
\newcommand{\bn}{\overline{\nabla}}
\newcommand{\bD}{\overline{D}}
\newcommand{\bT}{\overline{T}}
\newcommand{\tg}{\tilde{g}}		
\newcommand{\htg}{\hat{\tilde{g}}}
\newcommand{\SO}{\operatorname{SO}}
\newcommand{\Symp}{\operatorname{Sp}}
\newcommand{\SU}{\operatorname{SU}}
\newcommand{\Si}{\operatorname{Spin}}			
\newcommand{\Sr}{{\cal{S}}}		
\renewcommand{\phi}{\varphi}		

\newcommand{\hy}{{\Bbb H}^n}		
\newcommand{\hp}{{\Bbb H}_*^n}		
\newcommand{\B}{B^n}
\newcommand{\dM}{\partial M}
\renewcommand{\mod}{\operatorname{mod}}
\newcommand{\tR}{\widetilde{R}}
\newcommand{\Cf}{\frac{1}{2} (1-|x|^2)}
\newcommand{\half}{\frac{1}{2}}

\newcommand{\ihalf}{\frac{{\mathbf{i}}}{2}}
\newcommand{\bfi}{{\mathbf i}}
\newcommand{\quart}{\frac{1}{4}}

\newcommand{\hD}{\widehat{D}}
\newcommand{\hB}{\widehat{B}}
\newcommand{\hR}{\widehat{R}}
\newcommand{\hs}{\hat{s}}		
\newcommand{\hrho}{\hat{\rho}}
\newcommand{\hx}{\hat{x}}
\newcommand{\hf}{\hat{f}}

\newcommand{\tG}{\widetilde{\Gamma}}

\newcommand{\la}{\langle}	
\newcommand{\ra}{\rangle}	

\newcommand{\tr}{\operatorname{tr}}

\renewcommand{\qed}{\hfill~\begin{picture}(8,8)(0,0)
                \thicklines
                \put(0,.5){\line(1,0){8}}
                \put(7.5,0){\line(0,1){8}}
                \thinlines
                \put(0,8){\line(1,0){8}}
                \put(0,0){\line(1,0){8}}
                \put(8,0){\line(0,1){8}}
                \put(0,0){\line(0,1){8}}
                \end{picture}
                \bigskip
                }

\newcommand{\proof}{{\sc Proof:}}
\theoremstyle{plain}
\newtheorem{thm}{Theorem}[section]
\newtheorem{lemma}[thm]{Lemma}
\newtheorem{prop}[thm]{Proposition}
\newtheorem{cor}[thm]{Corollary}
\newtheorem{remark}[thm]{Remark}
\newtheorem{definition}[thm]{Definition}

\title[Asymptotically locally hyperbolic manifolds]{Scalar Curvature Rigidity for asymptotically locally hyperbolic manifolds}
\author[Lars Andersson]{Lars Andersson$^1$}
\thanks{Submitted to {\em Annals of Global Analysis and Geometry}}
\thanks{$^1$Research partially supported by the Swedish Natural Sciences Research Council,
contract no. M-AA/MA 04873-313 and the Crafoord foundation.}
\author[Mattias Dahl]{Mattias Dahl$^2$}
\thanks{$^2$Research partially supported by the Wallenberg foundation.}
\address{
Department of Mathematics \\
Royal Institute of Technology \\
S-100 44 Stockholm, Sweden\\
}
\email{larsa\char'100math.kth.se, dahl\char'100math.kth.se}


\begin{document}
\date{Aug. 14, 1996}
\begin{abstract}
Rigidity results for asymptotically locally hyperbolic manifolds with lower
bounds on scalar curvature are proved using spinor methods related to the
Witten proof of the positive mass theorem. The argument is based on a study
of the Dirac operator defined with respect to the Killing connection. 
The existence of asymptotic Killing spinors is related to the spin structure
on the end.
The expression for the mass is calculated and proven to vanish for
conformally compact Einstein manifolds with conformal boundary a spherical
space form, giving rigidity. 
In the 4-dimensional case, the signature of the manifold is related to the
spin structure on the end and explicit formulas for the relevant invariants
are given. 
\end{abstract}

\maketitle

{ \small
\tableofcontents
}

\section{Introduction}
\label{sec:intro}
Let $(M,g)$ be a complete spin manifold of dimension $n \geq 3$. We say that
$(M,g)$ is asymptotically locally hyperbolic (ALH) if $(M,g)$ has an end $E$ 
where the metric is asymptotic to $\hy/\Gamma$ where $\hy$ is hyperbolic
space and $\Gamma$ is a finite group acting by
isometries, see section \ref{sec:ALH}. In particular the sectional curvatures 
are asymptotic to $-1$ on $E$. In this paper we will prove some 
rigidity results for ALH manifolds whose scalar curvature $s$ satisfies 
$s \geq -n(n-1)$. The notion of an ALH manifold is analogous to the notion 
of an asymptotically locally Euclidean (ALE) manifold, a class of manifolds 
which has been extensively studied.

In \cite{min-oo:scalar} Min-Oo proved a scalar curvature rigidity result for
manifolds asymptotic to hyperbolic space. The argument was an adaptation of
the Witten proof of the positive mass theorem. The asymptotically hyperbolic
case is related to the proof of positivity of the Bondi mass. 

From the point of view of Riemannian geometry the positive mass theorem may
be viewed as a statement concerning asymptotically Euclidean 3-manifolds with
nonnegative scalar curvature. For
dimension higher than 3, we have the generalized positive action conjecture, 
for ALE manifolds with nonnegative scalar curvature, which was shown to be 
false by LeBrun \cite{lebrun:counter}. The positive mass argument fails when
the spin structure does not allow asymptotically parallel spinors.
Even given the existence of asymptotically parallell spinors the conclusion
drawn from the positive mass argument in the ALE case is in general just a
restriction on holonomy, see however \cite{kronheimer:torelli} for a complete
discussion of the $4$-dimensional case.

In the ALH case, it is natural to consider 
a modified connection which is such that its curvature
vanishes on hyperbolic space. The corresponding Dirac operator and
Lichnerowicz identities can then be used in the positive mass argument in
much the same way as in the Witten proof, the relevant mass being related to
the Bondi mass. 

There are different ways to modify the Levi-Civita connection so that 
it becomes adapted to hyperbolic space. In \cite{min-oo:scalar} a connection 
on an extended spinor bundle was studied.
This setup models the manifold as an umbilic hypersurface in a Minkowski
space of dimension $n+1$. Alternatively one may introduce a modified
connection $\hn$ defined by $\hn_X  \psi = \n_X \psi + \ihalf X \psi$ 
on the usual spinor bundle, which has as parallel spinors the
imaginary Killing spinors. This is the approach taken in the present paper. 

An important difference between Killing spinors and parallel spinors is that
the existence of an imaginary Killing spinor allows one to apply the powerful
structure theorem of H. Baum, see Theorem \ref{thm:baum} and in contrast to 
the ALE case, instead of restrictions on the holonomy we get rigidity.

If $M$ is ALH or ALE, then $M$ has an end $E$ diffeomorphic to 
$\Re \times (S^{n-1}/\Gamma)$. The existence of asymptotically parallel 
or asymptotically Killing spinors is a topological condition on the spin
structure on $E$.  When $M$ is 4-dimensional we use the Atiyah-Patodi-Singer 
index theorem to
relate the existence of asymptotic Killing spinors to the signature of $M$.
From this we derive conditions on the signature of $M$ from the ALH condition
and the lower bound on the scalar curvature. 

As a particular case we study conformally compact Einstein
manifolds. We calculate the expression for the mass for a conformally compact
manifold and in case the conformal boundary is covered by the round sphere we
prove that the mass vanishes, which leads to rigidity, see section
\ref{sec:conf:compact}.

\noindent{\bf Overview of this paper:}
This paper is organized as follows. In section \ref{sec:prelim} we discuss
spin structures on quotients and the problem of relating spinors 
defined with respect to different metrics. 
In section \ref{sec:Killing-main} we give the
necessary background on imaginary Killing spinors. In section 
\ref{sec:rigidity} we define the notion of ALH manifold and prove the basic
rigidity theorem, Theorem \ref{thm:main1}. The main step in the proof is Lemma 
\ref{lemma:min-oo} which is an adaptation of the argument in
\cite{min-oo:scalar} to the present situation.
In section \ref{sec:conf:compact} we apply these ideas to the particular case
of conformally compact Einstein manifolds and prove that when the conformal
boundary is a spherical spaceform, the mass vanishes which reduces the problem
of proving scalar curvature rigidity to a topological question of whether the
spin structure on the end admits an asymptotic Killing spinor. 
Finally, in section \ref{sec:4d} we discuss the relation between the topology
of $M$ and the spin structure on $E$ which leads up to Corollary
\ref{cor:4dALH} and Corollary \ref{cor:4dEin}.

\noindent{\bf Acknowledgements:} We thank R. Rubinstein for some enlightening
discussions and H. Baum and T. Kunstmann for helpful remarks on an 
earlier version. 

\section{Preliminaries}
\label{sec:prelim}

Let $(M,g)$ be an oriented Riemannian spin manifold of dimension $n \geq 3$ with 
spin structure $\Si(M,g)$ and let $\Sr(M,g)$
be the spinor bundle on $M$ associated to $\Si(M,g)$. 

\subsection{Spin structures on quotients}
\label{sec:spinquot}
Let $\Gamma$ be a group acting 
by orientation preserving isometries on $M$. 
An element $\gamma \in \Gamma$ acts on a frame $f$ by 
$f \mapsto \gamma_* f$. Assume that this action of $\Gamma$ on the frame
bundle lifts to an action on $\Si(M,g)$, that is we
have an action $s \mapsto \tilde \gamma s$ which projects to the action 
$f \mapsto \gamma_* f$.
Via the spin representation this defines an action on the spinor 
bundle where we denote the action of $\gamma$ by $\tilde \gamma$.

Assume that $\Gamma$ is a discrete group acting without fixed points. Then
$\Gamma$ has a lift if and only if $M/\Gamma$ is spin. In this case the
spin bundle on the quotient is given by 
$$
\Si(M/\Gamma) = \Si(M)/\Gamma
$$ 
and the associated spinor bundle is given by 
$$
\Sr(M/\Gamma) = \Sr(M)/\Gamma.
$$
This means that given a lift the sections of $\Sr(M/\Gamma)$ are precisely the
$\Gamma$-periodic sections of $\Sr(M)$.

Lifts of $\Gamma$ are classified by $\operatorname{Hom}(\Gamma,\Za_2)$ 
and in case $M$ is simply connected it follows from the isomorphism 
$H^1(M / \Gamma ; \Za_2) = \operatorname{Hom}(\Gamma,\Za_2)$ that this also 
classifies the spin structures on $M/\Gamma$.

\subsection{Comparing spinors for different metrics}
\label{comparespin}
Let $g,g'$ be Riemannian metrics on a  manifold $M$, and define the 
positive definite `gauge transformation' $A \in \operatorname{End}(TM)$ by
$$
g(AX,AY) = g'(X,Y),
$$
$$
g(AX,Y) = g(X,AY).
$$
Because of the first property $A$ will map ON-frames for $g'$ to ON-frames
for $g$, and thus $A$ induces a map $\SO(M,g') @>{A}>> \SO(M,g)$. If $M$ is spin
and we choose equivalent spin-structures for $g$ and $g'$ this can be lifted
to $\Si(M,g') @>{A}>> \Si(M,g) $. A spinor field for $g$ can be viewed as
a $\Si(n)$-equivariant map $\Si(M,g) @>{\phi}>> \Sr$, where $\Sr$ is the
spinor space, so the composition
$\phi \circ A$ is a map $\Si(M,g') @>{\phi}>> \Sr$ which also is 
$\Si(n)$-equivariant. This gives the extension of $A$
to a map $\Sr(M,g') @>{A}>> \Sr(M,g)$ which respects Clifford multiplication;
$$
A(X\cdot \phi) = (AX) \cdot (A\phi).
$$
Since the metric on the spinor bundle is given by a fixed Hermitean 
inner product on $\Sr$, $A$ defines a fibrewise isometry.
The above can be collected in a diagram.
$$
\begin{CD}
{\Si(M,g')} @>{A}>> {\Si(M,g)} @>{\phi}>> {\Sr} \\
   @VV{\pi}V @VV{\pi}V  \\
{\SO(M,g')} @>{A}>> {\SO(M,g)} \\
\end{CD}
$$

We will now look at the relation between the canonical covariant derivatives 
for $(M,g)$ and $(M,g')$.
Let $\n$ and $\n'$ be the Levi-Civita connections
for $g$ and $g'$, to be able to compare $\n$ and $\n'$ on the frame and 
spin bundles for $g$ we define a connection $\bn$ by
\begin{equation} \label{barnabla}
\bn X= A ( \n'  A^{-1}X).
\end{equation}
The connection $\bn$ is metric with respect to $g$ and has torsion
\begin{eqnarray}
\bT(X,Y) &=& \bn_X Y - \bn_Y X - [X,Y] \label{torsion} \\
\nonumber &=& -((\n'_X A)A^{-1}Y - (\n'_Y A)A^{-1}X).
\end{eqnarray}
Expressing the covariant derivative in 
terms of the Lie bracket and the metric we get 
\begin{equation} \label{difftorsion}
2g(\bn_X Y - \n_X Y,Z) = g(\bT(X,Y),Z) - g(\bT(X,Z),Y) - g(\bT(Y,Z),X).
\end{equation}
Next we compare $\n,\bn$ when lifted to the spinor bundle $\Sr(M,g)$. 
Let $\{ e_i\}$ be a local orthonormal frame for $g$, and 
let $\{ \sigma_{\alpha}\}$ be the corresponding local 
orthonormal frame of the spinor bundle.
Denote by $\omega_{ij},\overline{\omega}_{ij}$ the 
connection one-forms for $\n,\bn$ defined with respect to $\{e_i\}$,
$$
\omega_{ij} = g(\n e_i,e_j),
$$
$$
\overline{\omega}_{ij} = g(\bn e_i,e_j),
$$
then the covariant derivatives of $\phi = \phi^{\alpha} \sigma_{\alpha} $
are given by \cite[Thm 4.14]{LM}
$$
\n \phi = d\phi^{\alpha}\otimes \sigma_{\alpha} + 
\half \sum_{i<j} \omega_{ij} \otimes e_ie_j \phi,
$$
$$
\bn \phi =d\phi^{\alpha}\otimes \sigma_{\alpha} + 
\half \sum_{i<j} \overline{\omega}_{ij} \otimes e_ie_j \phi 
$$
and hence the difference between $\bn$ and $\n$ acting on $\phi$ is
\begin{equation}
\bn \phi - \n \phi = \half \sum_{i<j} (\overline{\omega}_{ij} -
\omega_{ij})\otimes e_ie_j \phi.
\end{equation}
Using (\ref{torsion}) and (\ref{difftorsion}) we can estimate
$$
|(\overline{\omega}_{ij}-\omega_{ij})(e_k)| \leq
C |A^{-1}| |\n'A|.
$$
We have proved the following lemma
\begin{lemma}\label{lemma:diff}
Let $Y$ be a vectorfield and let $\phi$ be a spinor (w.r.t the $g$ spin bundle),
then
\begin{equation}
\label{diffvector}
|\bn Y - \n Y| \leq C |A^{-1}| |\n'A| |Y|, 
\end{equation}
\begin{equation}
\label{diffnabla}
|\bn \phi - \n \phi| \leq C |A^{-1}| |\n'A| |\phi| 
\end{equation}
and
\begin{equation}
\label{diffdirac}
|\bD \phi - D \phi| \leq C |A^{-1}| |\n'A| |\phi|,
\end{equation}
where $D,\bD$ are the Dirac operators associated to the connections $\n, \bn$.
~\qed
\end{lemma}

\section{Killing spinors}
\label{sec:Killing-main}
\subsection{The Killing connection}
\label{sec:Killing}
Define a modified connection on the spinor bundle by
$$
\hn_X = \n_X + \ihalf X \cdot ,
$$
we call this the Killing connection, spinors parallel with respect to this connection are called imaginary Killing spinors. There is an analogous concept of
real Killing spinors. In this paper we will discuss only the imaginary
Killing spinors and write simply Killing spinors. 
The connection $\hn$ will respect the Clifford module structure of
$\Sr(M,g)$ if we define its action on the Clifford algebra bundle on $(M,g)$
so that 
$$
\hn_X Y = \n_X Y + \ihalf [X,Y],
$$
where $[ \, ,\, ]$ denotes the commutator in the Clifford algebra. 

The Dirac operator corresponding to $\hn$ is
$$
\hD 
= D - \ihalf n.
$$
The curvature of $\hn$ is 
\begin{eqnarray*}
\hR (X,Y)\phi &=& (\hn_X\hn_Y -\hn_Y\hn_X - \hn_{[X,Y]})\phi \\
&=& (\n_X\n_Y -\n_Y\n_X - \n_{[X,Y]})\phi \\
&&\hskip 0.2in + \ihalf ((\n_XY)-(\n_YX) -[X,Y])\phi -\frac{1}{4}(X\cdot Y - Y\cdot
X)\cdot
 \phi \\
&=& R(X,Y)\phi-\frac{1}{4}(X\cdot Y - Y\cdot X)\cdot\phi
\end{eqnarray*}
Through the natural isomorphism of Lie algebras 
${\frak{spin}}(n) \leftrightarrow {\frak{so}}(n)$ (\cite[\S I.6]{LM}) 
the second term corresponds to the skew-adjoint 
endomorphism $X \wedge Y$
of $TM$ defined by 
$$
(X \wedge Y)(Z) = \la X,Z \ra Y -\la Y,Z \ra X.
$$ 
This is also the curvature tensor $R_{-1}$ on $TM$ of constant sectional 
curvature $-1$, so we see that the curvature of the Killing connection
viewed as an ${\frak{so}}(n)$-valued two form is
\begin{equation}\label{eq:Rhat}
\hR(X,Y) = R(X,Y) - R_{-1}(X,Y).
\end{equation}

\subsection{Manifolds with Killing spinors}
This expression for the curvature of $\hn$ tells us that
if there is a local basis of Killing spinors then $\hR$ 
vanishes and $(M,g)$ is locally isometric to hyperbolic space.

If a manifold has one spinor parallel with respect to $\n$ it is Ricci-flat,
this follows from the formula 
$$
\sum_i e_i R (X,e_i)\phi = -\half \operatorname{Ric}(X)\cdot\phi.
$$
We calculate using (\ref{eq:Rhat}) for a Killing spinor $\phi$
$$
0 = \sum_i e_i \hR (X,e_i)\phi = -\half(\operatorname{Ric}(X) + (n-1)X
)\cdot\phi, 
$$
from which it follows that a manifold with a Killing spinor must be Einstein.
But for a complete manifold with a Killing spinor there is the following
result of H. Baum.
\begin{thm}[H. Baum 
\protect{\cite{baum:imaginary:killing}} 
\protect{\cite[Chapter 7]{baum:killing}}] 
\label{thm:baum}
A complete manifold $M$ has a Killing spi\-nor $\phi$ if and only if it is
isometric to $P \times \Re$ with the metric 
$ e^{2t}h+dt^2$ where $(P,h)$ is a
complete manifold with a parallel spinor.
~\qed
\end{thm} 

On hyperbolic space $\hy$ there is a full set of 
Killing spinors. We describe them using the ball model of $\hy$, 
which is $\B = \{ |x|<1 \}$ 
with the metric $g= \rho(x) ^{-2} g_0$ where $\rho(x) = \Cf$ and $g_0$ is the
standard flat metric. Since this is 
conformal to the flat metric on the ball we get an identification 
$\SO(\hy) = \B \times \SO(n)$ which gives trivializations 
$\Si(\hy) = \B \times \Si(n)$ and $\Sr(\hy) = \B \times \Sr$. 

\begin{prop}[\cite{baum:imaginary:killing} \cite{baum:killing}]\label{prop:hypkill}
In the above trivialization, the Killing spinors on $\hy$ are 
$$
\phi_u(x) = \rho(x)^{-\half} (1+\bfi x \cdot) u, 
$$ 
where $u \in \Sr$.
\end{prop}

\subsection{The Lichnerowicz formula}
Fix a spinorfield $\phi$ and define a one-form
$$ 
\hat \alpha(X) = \la (\hn_X +X \cdot \hD) \phi,\phi \ra. 
$$
A computation  gives
\begin{eqnarray*}
\operatorname{div} \hat \alpha 
&=&\frac{\hs}{4}|\phi |^2 + |\hn \phi|^2 -| \hD \phi|^2 .
\end{eqnarray*}
where $\hs = s + n(n-1)$ and $s$ is the scalar curvature of $(M,g)$.
Integrating over a manifold $M$ with boundary $\dM$ we get the
Lichnerowicz formula
\begin{equation} \label{lichn2}
\int_M ( \frac{\hs}{4}|\phi |^2 + |\hn \phi|^2 -| \hD \phi|^2) = \int_{\dM} \la (\hn_{\nu} +\nu \hD) \phi,\phi \ra,
\end{equation}
where $\nu$ is the outward normal of the boundary. 
We also note the integration by parts
formula for the Dirac operator $\hD$.
\begin{equation} \label{intpart2}
\int_M \la (\hD + \bfi n) \phi,\psi \ra = \int_{\dM} \la \nu \phi, \psi \ra
+ \int_M \la \phi,\hD\psi \ra.
\end{equation}

\subsection{A little analysis}
We end this section by proving some analytical results using the
formulas derived above. Let $M$ be a complete spin manifold with 
$\hs$ non-negative and bounded. 
Let $C^{\infty}_0\Sr(M)$ be the space of smooth sections with compact support
and define a sesquilinear form 
$$
B(\phi,\psi) = \int_M \la \hD \phi,\hD \psi \ra.
$$
which is obviously bounded. $B$ is Hermitean and using (\ref{lichn2}) we get 
\begin{eqnarray*}
B(\phi,\phi) &=& \int_M\frac{\hs}{4} |\phi|^2 + |\hn\phi|^2 \\
&=& \int_M\frac{\hs +n}{4}|\phi|^2 + |\n \phi|^2 
\end{eqnarray*}
Define a scalar product $(\cdot ,\cdot)_1$  by 
$$
(\phi,\psi)_1 = \int_M \la \n \phi,\n \psi \ra  + \frac{n}{4}\la \phi,\psi \ra  
$$
and let $H^1\Sr(M)$ be the closure of $C^{\infty}_0 \Sr(M)$ with respect 
to the norm 
$|| \phi||_1 = (\phi, \phi)_1^{1/2}$. 
Then $(H^1\Sr(M), (\cdot,\cdot)_1)$ is a Hilbert space. 
If $\hs$ is bounded, $B$ extends to $H^1\Sr(M)$ as a bounded sesquilinear form 
and if $\hs$ is non-negative we see that $B$ is coercive on $H^1\Sr(M)$. 

We have shown that $B$ satisfies the conditions of the Theorem of Lax-Milgram,
the conclusion is 
\begin{prop} \label{lax}
Let $(M,g)$ be as above and assume that $\hs$ is nonnegative and bounded.
Then for every bounded linear functional $l$ on $H^1\Sr(M)$ there is a unique 
$\phi\in H^1\Sr(M)$ so that for all $\psi \in H^1\Sr(M)$ 
$$
B(\phi,\psi) = l(\psi).
$$
~\qed
\end{prop}  
\begin{remark}
In the rest of the paper we will assume that $\hs$ is bounded. 
This assumption is technical and can be removed by
using appropriately weighted spaces of sections. 
~\qed
\end{remark}
\section{Rigidity for ALH manifolds}
\label{sec:rigidity}
In this section we consider manifolds satisfying a condition which is
analogous to the vanishing of the Bondi mass. By the positive mass argument
this leads to rigidity theorems. 

\subsection{Asymptotic Killing spinors}

Let $(M,g)$ be a complete Riemannian manifold with one end $E$ and a
diffeomorphism $N \times (0,\infty) @>{\phi}>> E$ such that 
$\phi^* g = L^2 dR^2 + h$
where $h$ is the induced metric on $N_R = \phi(N \times \{ R\})$ and $L > 0$. 
Let $L_{\text{min}} (t) = \min_{N_t} L$. We say that $\phi$ gives a nondegenerate
foliation at infinity if $\int_0^{\infty} L_{\text{min}}(t) dt = \infty$. Let
$\{M_R\}$ be a family of compact manifolds with boundary exhausting $M$. We
say that that $\{M_R\}$ is nondegenerate if there is a nondegenerate
foliation at infinity such that $\partial M_R = N_R$. 

We have
\begin{prop}\label{prop:tobias} Let $(M,g)$ be a complete Riemannian manifold 
with a nondegenerate exhausting family $M_R$. Let $f$ be a nonnegative 
integrable founction on $M$. Then 
$\liminf_{R \to \infty} \int_{\dM_R} f = 0$.
~\qed
\end{prop}

The generalization to
manifolds with more than one end is trivial.

\begin{definition} \label{def:mass}
Let $(M,g)$ be a complete Riemannian spin manifold.
\begin{enumerate}
\item 
We say that a spinor $\phi_0$ is an 
asymptotic Killing spinor on $M$ if $\hn\phi_0 \in L^2(T^{*}\otimes\Sr)(M)$, 
$\phi_0 \notin H^1\Sr(M)$. 
\item 
If $\phi_0$ is an asymptotic Killing spinor we say
that the mass of $M$ (w.r.t. $\phi_0$) is zero if there is a nondegenerate 
family $\{M_R\}$ exhausting $M$ 
so that either 
\begin{equation}\label{eq:masszerolim}
\lim_{R \to \infty} \int_{\dM_R} \la (\hn_{\nu} +\nu \hD) \phi_0,\phi_0 \ra = 0
\end{equation}
or 
\begin{equation}\label{eq:masszeroL1}
\la (\hn_{\nu} +\nu \hD) \phi_0,\phi_0 \ra \in L^1(M).
\end{equation}
\end{enumerate}
~\qed
\end{definition}

\begin{remark} In the setting of Definition \ref{def:mass}, if the limit
$$
\lim_{R \to \infty} \int_{\dM_R} \la (\hn_{\nu} +\nu \hD) \phi_0,\phi_0 \ra 
$$
exists, then it is a natural to think of this as an analog of a component of
the Bondi momentum in general relativity.
~\qed
\end{remark}

The following Lemma is the essential step in the positive mass argument.
\begin{lemma}
\label{lemma:basic}
If $M$ has an asymptotic Killing spinor,
mass zero in the sense of Definition \ref{def:mass} and
$\hs \geq 0$, then $M$ has a Killing spinor.
\end{lemma}
\proof Since $\hD\phi_0 \in L^2\Sr(M)$ the linear functional
$$
l(\psi) = \int_M \la \hD\phi_0,\hD\psi\ra
$$
is bounded on $H^1\Sr(M)$, so by Proposition \ref{lax} there is a unique 
$\phi_1 \in H^1\Sr(M)$ so that for all $\psi \in H^1\Sr(M)$ 
$$
\int_M \la \hD \phi_1,\hD \psi \ra = B(\phi,\psi) = l(\psi) = 
\int_M \la \hD\phi_0,\hD\psi\ra.
$$
Set $\phi = \phi_0 - \phi_1$, then  
$$
\int_M \la \hD \phi,\hD \psi \ra = 0,
$$
we will show that this implies $\hD\phi =0$. Set $a=\hD\phi$. 
Integrating by parts using (\ref{intpart2}) we get for any 
$\psi \in C^{\infty}_0\Sr(M)$
$$
0  = \int_M \la a,\hD \psi \ra  = \int_M \la (\hD+\bfi n)a,\psi \ra, 
$$
so 
\begin{equation}
\label{D+in}
(\hD+\bfi n)a = 0  .
\end{equation}
By elliptic regularity we get that $a$ is smooth and by (\ref{D+in})
$\hD^k a \in L^2\Sr(M)$ for any integer $k \geq 0$. 
Therefore the usual cutoff argument shows that
\begin{eqnarray*}
\int_M \la \hD a,\hD a\ra &=& \int_M \la(\hD+\bfi n)\hD a,a\ra\\
&=& \int_M \la\hD(\hD+\bfi n)a,a\ra\\
&=& 0
\end{eqnarray*}
and thus $\hD a =0$ which together with $(\hD+\bfi n)a =0$ implies 
$a=\hD\phi =0$.
Let $M_R$ be as in definition \ref{def:mass}. We use (\ref{lichn2}) 
to get 
\begin{eqnarray}
\int_{M_R} |\hn \phi|^2 &\leq& \int_{M_R} ( \frac{\hs}{4}|\phi |^2 
+ |\hn \phi|^2) \nonumber \\
&=&
\int_{\partial M_R} \la (\hn_{\nu} +\nu \hD) \phi,\phi \ra. \label{boundterm}
\end{eqnarray}
Write $\hB = \hn_{\nu} + \nu \hD $. 
Since $\hB$ is self adjoint on 
$\partial M_R$ we may write 
(\ref{boundterm}) in the form 
\begin{equation}\label{eq:bdryterm}
\int_{\partial M_R} \la \hB \phi_0,\phi_0 \ra 
+ \int_{\partial M_R} \la \hB \phi_1,\phi_1 \ra
+ \int_{\partial M_R} \la \hB \phi_0 , \phi_1 \ra + 
\int_{\partial M_R} \la \phi_1 , \hB \phi_0 \ra .
\end{equation}
In case (\ref{eq:masszerolim}) holds, 
the first term tends to zero as $R \to \infty$ from the vanishing
of the mass w.r.t. $\phi_0$. By assumption, $\hn \phi_0 \in 
L^2(T^{*}\otimes\Sr)(M)$ and by
construction $\phi_1 \in H^1\Sr(M)$, therefore we have $|\la \hB \phi_0,
\phi_1 \ra | \in L^1(M)$ and $|\la \hB \phi_1 , \phi_1 \ra | \in L^1(M)$ 
and by Proposition \ref{prop:tobias}  there is a sequence $\{R_j\}$ such that 
the last three terms tend to zero as $j \to \infty$. 
This tells us that $\hn\phi =0$ and since $\phi_0 \notin H^1\Sr(M)$
by assumption we have $\phi \neq 0$. Finally in case (\ref{eq:masszeroL1})
holds, we may apply Proposition \ref{prop:tobias} as above to all the terms
in (\ref{eq:bdryterm}).
\qed

\subsection{Asymptotically locally hyperbolic manifolds}
\label{sec:ALH}
\begin{definition}\label{def:SA}
Two metrics $g$ and $g'$ are strongly asymptotic if the 
gauge--trans\-for\-mation $A$
of Section \ref{comparespin} satisfies
\begin{enumerate}\closeup
\item \label{point:SA1} 
There is a $k$ so that for all $X \in TM, |X| =1$ we have $k^{-1} 
\leq |AX| \leq k$, 
\item \label{point:SAnew} $ (|\n' A|^2 + |A-Id|^2)^{\half} \in L^2(M) \cap
L^1(M)$
w.r.t. the measure $e^rd\operatorname{vol}(g)$,
\item \label{point:SA4} 
$ |{\n'}^2 A| \in L^1(M)$ 
\mnote{\lars changed point labeled SA4 (point 3) in definition \ref{def:SA}}
\mnote{ $r$ or $r'$!!}
\end{enumerate} 
where $r$ is the $g'$-distance from a fixed point.
~\qed
\end{definition}
In the case when the model metric $g'$ is hyperbolic, this leads to
the concept of asymptotically locally hyperbolic manifolds,
analogous to the notion of asymptotically locally Euclidean spaces.
\begin{definition}\label{def:ALH}
Let $\Gamma$ be a finite group acting freely and linearly on the sphere and
let $\hp$ be $\hy\setminus B$ where $B$ is a closed ball centered at the
origin.
Suppose that $(M,g)$ has an end $E$ for which there is an 
identification of with $\hp/\Gamma$
where $\Gamma$ acts in the natural
way on the ball model of hyperbolic space.  
Let $g'$ be a metric on $M$ which on the end $E$ is isometric to 
$\hp/\Gamma$ with the hyperbolic metric.
Then we say that $(M,g)$ is asymptotically locally hyperbolic with the asymptotically locally hyperbolic end $E$
with group $\Gamma$ if $g$ and $g'$ 
are strongly asymptotic on $E$.
In case $\Gamma = \{ 1\}$ we call $M$ asymptotically hyperbolic.
~\qed
\end{definition}
As in the above definition let $\Gamma$ be a finite
subgroup of $\SO(n)$ acting freely on the sphere. We will now consider the
conditions for existence of Killing spinors on the quotient 
$\hp /\Gamma$. In the trivialization 
used in Proposition \ref{prop:hypkill} $\gamma \in \Gamma$ acts on 
$\SO(\hp)$ as
$$
(x,f) @>{\gamma}>> (\gamma(x),\gamma f).
$$
Assume that $\hp/\Gamma$ is spin. Then by the discussion in Section
\ref{sec:spinquot} there is a bijective lift of $\Gamma$ to a subgroup $\tG
\subset \Si(n)$,
which specifies the action of $\Gamma$ on the spin and spinor bundles of $\hp$.
The sections of $\Sr(\hp /\Gamma)$ are naturally identified with
the $\Gamma$-periodic sections of $\Sr(\hp)$.
Since the Killing condition is local the Killing spinors on 
$\hp /\Gamma$ are precisely given by the spinors $\phi_u$ in 
Proposition \ref{prop:hypkill} which are also $\Gamma$-periodic.    
The $\phi_u$ are acted on as follows
$$
(x, \phi_u(x) ) @>{\gamma}>> 
(\gamma(x), \rho(\gamma(x))^{-\half} \tilde{\gamma} (1+\bfi x \cdot) u )
$$
and 
$$ 
\tilde{\gamma} (1+\bfi x \cdot) u = (1+ \bfi  \tilde{\gamma} \cdot x \cdot 
\tilde{\gamma}^{-1} \cdot) \tilde{\gamma} \cdot u.
$$
If we view $\Si(n)$ as a subgroup of the Clifford algebra the
conjugation by $\tilde{\gamma}$ is just the projection to $\gamma$ in $\SO(n)$.
Therefore, if we assume that the fixed element $u$ in $\Sr$ satisfies
$$ 
\tilde{\gamma} \cdot u = u 
$$
for all $\tilde{\gamma}$ the Killing spinor will be $\Gamma$-periodic
$$ 
(x, \phi_u(x) ) @>{\gamma}>> 
(\gamma(x), \rho(\gamma(x))^{-\half} (1+\bfi \gamma(x) \cdot) u )
= (\gamma(x), \phi_u(\gamma(x))).
$$
Conversely if $\phi_u$ is $\Gamma$-periodic we must have 
$\tilde{\gamma} \cdot u = u$ for all $\tilde{\gamma}$.
The Killing spinors on the quotient thus correspond precisely to the spinors 
$u \in \Sr$ which 
are fixed by the spin-representation of the lifted group $\tG$, they depend 
both on $\Gamma$ and via the choice of lift the spin structure on the quotient.
We have proved the following Proposition.
\begin{prop}
\label{prop:quotient1}
Let $\Gamma$ be a subgroup of $\SO(n)$ acting freely on the sphere and let
$\tG$ be a bijective lift to $\Si(n)$, then the Killing spinors on
$\hp /\Gamma$ with the spin-structure defined by $\tG$ correspond to the 
$\phi_u$ given by Proposition \ref{prop:hypkill} for spinors $u \in \Sr$ 
which are fixed by the spinor representation of $\tG$.
~\qed
\end{prop}
Now let $M$ be an asymptotically locally hyperbolic spin manifold with
group $\Gamma$ and fix a spin structure on $M$. The restriction of the 
spin structure to the locally hyperbolic end will be 
equivalent to the spin structure on $\hy / \Gamma$ defined by some lift 
$\tG$ of $\Gamma$. Thus the number of asymptotic Killing spinors on $M$
is via the lift $\tG$ controlled by spin structure, a purely topological
condition.     

\subsection{A rigidity theorem}
We have the following rigidity theorem.
\begin{thm}
\label{thm:main1}
Let $M$ complete spin manifold with an asymptotically locally hyperbolic 
end $E$ with group $\Gamma$. 
Suppose that the spin structure on $E$ is equivalent to the spin structure 
on $\hp /\Gamma$ defined by the lift $\tG$ of $\Gamma$. 
If $\tG$ fixes some non--zero spinor $u \in \Sr$ and if 
$\hs \geq 0$ then $ \Gamma = \{ 1 \}$ and $M$ is isometric to $\hy$.
\end{thm}
\begin{remark} The condition on  the spin structure can be formulated as the
vanishing of a relative Stiefel--Whitney class, see
\cite{andersson:dahl:future}.
~\qed
\end{remark}
We begin by proving the existence of an asymptotic Killing spinor with mass
zero. 
\begin{lemma}
\label{lemma:min-oo}
Let $M$ complete spin manifold with an asymptotically locally hyperbolic
end $E$ with group $\Gamma$.
Suppose that the spin structure on $E$ is equivalent to the spin structure
on $\hp /\Gamma$ defined by the lift $\tG$ of $\Gamma$ and 
$\tG$ fixes some non--zero spinor $u \in \Sr$. 
Then $M$ has an asymptotic Killing spinor with mass zero.
\end{lemma}
\proof 
Let $E$ be the locally hyperbolic end and let $g'$ be the hyperbolic metric. 
On the end $\hp /\Gamma$ with the hyperbolic metric there is by
Proposition \ref{prop:quotient1} a Killing spinor $\phi'$. 
Since the spin-structures match up $\phi'$ can be moved to a spinor 
$A\phi'$  on $E$ where $A$ is the gauge transformation of section 
\ref{comparespin}. 

Let $f$ be a smooth function with $\operatorname{supp}(df)$ compact, $f=0$ 
outside $E$ and $f=1$ at infinity. 
Define the spinor $\phi_0 = fA\phi' $ on $(M,g)$. We are going to show 
that $\phi_0$ is an asymptotic Killing spinor. 

First  we show that $\hn \phi_0 \in L^2\Sr(M)$.
\begin{eqnarray} \label{nphi}
\hn_X \phi_0 &=& (\n_X + \ihalf X ) f A \phi' \\
\nonumber &=& (Xf) A\phi' + f ( \n_X  A\phi' + \ihalf X A \phi') \\
\nonumber &=& (Xf)A\phi' + f(\n_X - \bn_X ) A\phi'
+ f ( \bn_X + \ihalf X)A\phi'
\end{eqnarray}
where the first term is compactly supported and the second can be estimated
using (\ref{diffdirac}), since $\hn'\phi' =0$ the third term is $f$ times
\begin{eqnarray*}
(\bn_X + \ihalf X ) A \phi' &=& A \n'_X \phi' +
\ihalf X A \phi' \\
&=& -\ihalf ( AX - X) A\phi' \\
&=& - \ihalf (A - Id)XA\phi'.
\end{eqnarray*}
Thus we get for large $r$ 
$$
|\hn \phi_0 | 
\leq C(|\n'A|^2 + |A - Id|^2)|\phi_0|^2.
$$
Next we estimate $|\phi_0|^2$, since $\hn'\phi' = 0$ we have for $X$ 
with $|X|' =1$
$$
|X({|\phi'|}^2)| \leq 2|\la \n'_X\phi',\phi' \ra | = 
|\la X\cdot \phi',\phi' \ra | \leq {|\phi'|}^2, 
$$
if we take $X = \frac{\partial}{\partial r}$ and integrate along a radial
geodesic this gives
$$
C^{-1}e^{-r} \leq {|\phi'|}^2 \leq Ce^r
$$
and since $A$ is an isometry on spinors, we have asymptotically that 
\begin{equation} \label{normgrowth}
C^{-1}e^{-r} \leq |\phi_0|^2 \leq Ce^r.
\end{equation}
Thus we have proved 
\begin{eqnarray} 
\label{normgrowth1}
 |\hn_X \phi_0|^2 &\leq& C(|\n' A|^2 + |A-Id|^2)e^r |X|^2, \\
\label{normgrowth2}
 |\hD\phi_0|^2 &\leq& C(|\n' A|^2 + |A-Id|^2)e^r, 
\end{eqnarray}
which together with point \ref{point:SAnew} of Definition \ref{def:SA} tells us that 
$\hn\phi_0$ is in $L^2(T^{*}\otimes\Sr)(M)$.
Next we have to show that the mass is zero. Let $R$ be large and let 
$S_R = \{x \in E; r(x) = R\} $ denote the distance sphere on $E$ with outward
normal $\nu$.

From (\ref{normgrowth1}),(\ref{normgrowth2}) we get 
\begin{eqnarray*}
|\la (\hn_{\nu} +\nu \hD) \phi_0,\phi_0 \ra| &\leq&
C (|\n' A|^2 + |A-Id|^2)^{\half}|A\phi'|^2 \\ &\leq&
C (|\n' A|^2 + |A-Id|^2)^{\half}e^r
\end{eqnarray*}
and by point \ref{point:SAnew} of Definition \ref{def:SA} this is in $L^1(M)$. 

It remains to check that 
$\phi_0 \notin H^1\Sr(M)$. From (\ref{normgrowth}) it follows that 
$|\phi_0|^2 \geq C e^{-r}$ for large $r$ and hence $\phi_0 \notin L^2\Sr(M)$, 
since the sphere $S_R$ has area of the order $e^{(n-1)R}$.

This completes the proof of Lemma \ref{lemma:min-oo}.
~\qed

We have now showed that $\phi_0$ is an asymptotic Killing spinor with mass
zero and by 
Lemma \ref{lemma:basic} we conclude that 
$M$ has a Killing spinor $\phi$ and therefore $(M,g)$ is of the form
described in Theorem \ref{thm:baum} $(M,g)$.

 
The proof of Theorem \ref{thm:main1} is concluded by the following Lemma.
\begin{lemma}
\label{lemma:done} 
Let $(M,g)$ be a warped product $(P \times \Re, e^{2t} h + dt^2)$ where 
$(P,h)$ is a complete Riemannian manifold. 
Suppose that $M$ has an asymptotically locally hyperbolic end. 
Then $M$ is isometric to $\hy$.
\end{lemma}
\proof 
Let $\xi,\eta$ be orthonormal vectors
orthogonal to $\frac{\partial}{\partial t}$, then $e^{t}\xi,e^{t}\eta$
are orthonormal w.r.t. $h$. The sectional curvature on the plane spanned by 
$\xi,\eta$ is (\cite[Chapter 7, Prop. 42]{oneill})
$$
K(\xi,\eta) = e^{-2t}K_0(\xi,\eta) - 1,
$$
where $K_0$ is the sectional curvature of $h$. On the other hand we 
know that there is an end on which $g(AX,AY)$ equals the hyperbolic metric
$g'(X,Y)$ of constant sectional curvature $-1$ where $A$ is as in 
Definition \ref{def:SA}. 
The curvature tensor of the hyperbolic metric $g'$ satisfies 
$$
R'(X,Y)Z = g'(X,Z)Y - g'(Y,Z)X.
$$
From (\ref{barnabla}) we see that $R'$ is related to 
the curvature $\overline{R}$ of $\bn$ by
$$
\overline{R}(X,Y)Z = A \circ R'(X,Y) \circ A^{-1}Z.
$$
It follows that 
$$
\overline{K}(\xi,\eta) = g(\xi, A \eta)^2 -g(\xi,A\xi)g(\eta,A\eta).
$$
and
$$
|\overline{K}(\xi,\eta) +1| \leq C |A-Id| .
$$
Define the tensor $\delta$ by
$$
\bn_X = \n_X + \delta_X .
$$
The curvature tensors for $\bn$ and $\n$ are related by
$$
\overline{R}(X,Y)Z = R(X,Y)Z 
+ (\bn_X\delta_Y - \bn_Y\delta_X - [\delta_x,\delta_y] + \delta_{[X,Y]})Z .
$$
Using (\ref{torsion}),(\ref{difftorsion}) and point \ref{point:SA1} of Definition
\ref{def:SA} we get 
$$
|\overline{K}(\xi,\eta)-K(\xi,\eta)| \leq C(|\n' A|^2 + |\n' A| + |{\n'}^2 A|)
$$
and 
\begin{eqnarray*}
e^{-2t}|K_0| &=& |K + 1 - \overline{K} + \overline{K}|\\
&\leq& |K- \overline{K}|+|\overline{K}+1|\\
&\leq& C(|{\n'}^2A| + |\n' A|^2 + |\n'A| + |A-Id|).
\end{eqnarray*}
For $x \in P$, let $K_{max}(x)$ be the maximum of $|K_0|$ over the two--planes
at $x$. Then we get by
integrating over the slice $t=T$ 
\begin{eqnarray*}
&& e^{-2T}e^{(n-1)T}\int_{P} |K_{max}| d\operatorname{vol}(h)
\\
&& \hskip .3in 
\leq C \int_{t=T} \left ( |{\n'}^2A| + |\n' A|^2 + |\n' A| + |A-Id|  \right )  
\frac{\partial}{\partial t}  \lrcorner d\operatorname{vol}(g) .
\end{eqnarray*}
By points  \ref{point:SAnew} and \ref{point:SA4} of Definition \ref{def:SA} 
\mnote{\lars added ref to point \ref{point:SA4} of Definition \ref{def:SA}}
and Proposition
\ref{prop:tobias}, there is a sequence $\{T_j\}$ tending to $\infty$ so
that the right hand side vanishes as $j \to \infty$. Therefore since the
integral on the left hand side is independent of $T$ and $n \geq 3$
it must be equal to zero. 

This shows that $h$ is a complete flat metric and since the fundamental 
group of the end is finite 
we must have that $(P,h)$ is $\Re^n$ with the flat metric and hence $M=\hy$. 
~\qed 

\section{Rigidity for conformally compact Einstein manifolds}
\label{sec:conf:compact}
In this section we prove that the mass for conformally compact Einstein
manifolds with conformal boundary a spherical space form vanishes. This means
we can apply the results of section \ref{sec:rigidity}.

\subsection{Conformally compact Einstein manifolds}
\begin{definition} \label{def:cc}
 A Riemannian manifold $(M,g)$ is called conformally compact 
if $M$ is the interior of a compact manifold $\widetilde{M}$ with 
boundary $\partial M$ and $g$ is conformal to a smooth metric $\tg$ on
$\widetilde{M}$,
$$
g = \rho ^{-2} \tg,
$$
where $\rho$ is a defining function for the boundary, that is 
$\dM = \{ x : \rho(x) =0 \}$ and $d\rho \neq 0$ on $\dM$.
\end{definition}
In the work of Fefferman and Graham \cite{fefferman:graham} conformally
compact Einstein metrics which are such that the conformal background metric
$\tg$ is an even function of the defining function are studied near $\partial M$.
They prove that in case $M$ is of even dimension, the formal power series
expansion of $\tg$ in terms of $\rho$ is uniquely determined, while in the
odd dimensional case, there is an undetermined term at the $n-1$:st order.
Here we study the case when the conformal boundary is locally conformally
flat and covered by a sphere. Using a combination of the spinor arguments
above and an analysis of the power series expansion of $\tg$ we prove the
following rigidity theorem.
\mnote{compare to infinitesimal rigidity of Graham and Lee}
\mnote{only need one component of boundary}
\begin{thm}\label{thm:conf}
Let $(M,g)$ be a conformally compact Einstein spin manifold of dimension $n$
for which the boundary of the conformal background metric is isometric
to $S^{n-1}$, then $(M,g)$ is isometric to $\hy$.
\end{thm}
\mnote{discuss dimension 2 and 3}
Let $(S^{n-1}, h_0)$ be the standard sphere of dimension $n-1$. Hyperbolic
space $\hy$ can be represented as $\Re_+ \times S^{n-1}$ with the metric
$$
\sinh^{-2}(x)( dx^2 + h_0)  .
$$
We will use a conformal gauge change $\tg \to \theta^2 \tg$ and
$\rho \to \theta\rho$ where $\theta$ is a positive function to 
put a general conformally compact Einstein metric on a form similar 
to this.  

\begin{lemma}\label{lemma:confgauge}
Let $(M,g)$ be as in Theorem \ref{thm:conf}.
There is a conformal gauge
change on a collar neighborhood $U$ of  $\dM$ so that $\rho = \sinh(x)$ where 
$x$ is the distance to $\partial M$ in $\tg$.
\end{lemma}

\proof Given are $M,\rho,\tg$ such that $g=\rho^{-2}\tg$ is Einstein.
The relation between the Ricci-tensors of $g$ and $\tg$ is
\begin{equation}\label{ricciconftrans}
\Ric = \tRic + \rho^{-1} ( (n-2) \tn d\rho + \tr_{\tg} ( \tn d\rho) \tg)
- (n-1) \rho^{-2} \tg(d\rho, d\rho) \tg.
\end{equation}
Since $\Ric = -(n-1)g$ this gives 
$$
-(n-1)\tg = \rho^{2}\tRic + \rho ( (n-2) \tn d\rho + 
\tr_{\tg} ( \tn d\rho) \tg) - (n-1) \tg(d\rho, d\rho) \tg,
$$
which implies that $\tg(d\rho, d\rho)=1$ on $\dM$ so 
$1-\tg(d\rho, d\rho) = \rho a$ where $a$ is some smooth function. 
We are going to find $\theta$ so that $\hrho := \theta\rho = \sinh(\hx)$
where $\hx$ is the distance to the boundary in the metric $\htg:=\theta^2\tg$,
or equivalently so that $\hf:=\sinh^{-1}(\hrho)=\hx$. We begin by solving
for a $\theta$ such that
$$
|d\hf|_{\htg} = 1.
$$
Written in terms of the given data this means that $\theta$ must 
solve the equation
$$
\rho^2\tg(d\theta,d\theta) + 2 \rho\theta\tg(d\theta,d\rho) =
\theta^4\rho^2 + \theta^2(1-\tg(d\rho,d\rho)),
$$
or 
$$
\rho\tg(d\theta,d\theta) + 2 \theta\tg(d\theta,d\rho) =
\theta^4\rho + \theta^2 a.
$$ 
One verifies easily that this equation with the initial data
$\theta=1$ on $\dM$ satisfies the conditions for existence 
of a solution in a neighbourhood of $\dM$ \cite[volume 5,  pp. 39-40]{spivak}.
In a small enough neighbourhood $U$ of $\dM$ the solution is
positive and we continue this to a positive
function $\theta$ on all of $M$ for which $|d\hf|_{\htg} = 1$ on $U$. 
This implies that the gradient curves of $\hf$ on $U$ are unit-speed
geodesics with respect to $\htg$ and since $\hf=0$ precisely 
on $\dM$ we have $\hf = \hx$ on $U$. 
\qed 

Let $\rho,x,\tg$ be as in Lemma \ref{lemma:confgauge} and let $\dM_x$ be the 
level surfaces of $x$. Then
$\tg$ is of the form $\tg = dx^2 + h$ where $h$ is the induced metric on
$\dM_x$. The normal vector field to the foliation $\dM_x$ w.r.t. $\tg$ is
$\eta = \frac{\partial}{\partial x} = \operatorname{grad}(x)$
and the second fundamental form of $\dM_x$ is
$$
\lambda(X,Y) = \tg(\eta, \tn_X Y) = -\half ( {\cal{L}}_{\eta} \tg ) (X,Y) 
$$
Observe that this identity makes sense for $X,Y$ not tangential
to $\dM_x$ by defining $\lambda(\eta, \cdot) =0$.

We define a hyperbolic background metric $g'$ on $U$ by setting
$$
g' = \sinh^{-2}(x)\tg'
$$
where $\tg'$ is the cylindrical metric $dx^2 + h_0$.
There is a Killing spinor $\phi'$ with respect to $g'$ 
defined on $U$. Let $A$ be the gauge transformation between 
$g$ and $g'$ and let $\phi_0 = fA\phi'$ where $f$ is 
a cut-off function. 
The proof of Theorem \ref{thm:conf} given in section \ref{sec:proof}
will consist of showing that $\phi_0$ is an asymptotic Killing spinor and 
that the mass of $(M,g)$ with respect to $\phi_0$ is zero. 
We begin by computing the asymptotics of the metric $\tg$.  

\begin{lemma}\label{lemma:A}
Let $n > 3$.
The gauge transformation $A$ satisfies 
\begin{eqnarray*}
(A - Id ) \eta &=& 0, \\
\tr_{h_0}(A - Id ) &=& O(x^{2n-2}), \\
A - Id  &=& O(x^{n-1}), \\
| \n A |_{g} &=& O(x^{n-1}), \\
| \n^2 A |_{g} &=& O(x^{n-1}). 
\end{eqnarray*}
Further, $x^{-(n-1)}(A - Id)$ is smooth up to $\partial M$. 
\end{lemma}
\proof
By Taylor's theorem we may expand $h$ in a asymptotic series
\begin{equation} \label{texp}
h = \sum_{k=0}^{\infty} \frac{x^k}{k!}h^{(k)}.
\end{equation}
where $h^{(0)} = h_0$ and $h^{(k)}$ are symmetric tensors on $\partial M$.
Then
\begin{equation} \label{second}
\lambda(X,Y)  = -\half \sum_{k=0}^{\infty} \frac{x^k}{k!}h^{(k+1)} (X,Y).
\end{equation}
The Hessian of $x$ equals $-\lambda$ so 
$$
\tn d\rho = \rho''(x) dx^2 - \rho' \lambda
$$
and, by assumption, the metric $g$ satisfies $\Ric =  -(n-1)g$. 
Using this and $\rho =\sinh(x)$ in (\ref{ricciconftrans}) gives us
the identity
\begin{equation}\label{confein}
\rho( \tRic + (n-2) dx^2  - (n-2) \tg ) = \rho' ( (n-2)  \lambda + \tr_{\tg}
\lambda  \tg ).
\end{equation}
Inserting (\ref{second}) in (\ref{confein}) and setting $x=0$
we get $h^{(1)} =0$. 

Recall the second variation formula 
\mnote{ref for this?} 
\begin{equation} \label{lieetal}
({\cal{L}}_{\eta} \lambda)(X,Y) = - \lambda \times \lambda (X,Y) 
+ {\cal{R}}_{\eta}(X,Y),
\end{equation}
where ${\cal{R}}_{\eta}(X,Y) = \tg(\tR(X,\eta),\eta,Y)$ and $\times$  
denotes the product  on 2-tensors defined by
$
a \times b (X,Y) = \sum_i a(X,\te_i)b(\te_i,Y)
$
where $\{ \te_i\}$ is an ON frame w.r.t. $\tg$.
Now let $\ric$ denote the Ricci tensor of the induced metric $h$ on $\dM_x$, 
then contracting the Gauss equation gives
\begin{equation} \label{Ricric}
\tRic = \ric + {\cal{R}}_{\eta} + \lambda \times \lambda - \tr_{\tg}(\lambda) \lambda
\end{equation}
on $\dM_x$.
Finally we plug (\ref{lieetal}) and (\ref{Ricric}) into equation 
(\ref{confein}) to get 
\begin{equation} \label{main}
\rho(\ric - (n-2) h + {\cal{L}}_{\eta} \lambda + 2 \lambda \times \lambda 
- \tr_{\tg}(\lambda) \lambda)
= \rho'( (n-2) \lambda + \tr_{\tg}(\lambda) h ).
\end{equation}
Further, the normal part of equation (\ref{confein}) is
$$
\rho \tRic ( \eta,\eta) = \rho'( (n-2) \lambda(\eta,\eta) 
+ \tr_{\tg} (\lambda)\tg(\eta,\eta))  = \rho' \tr_{\tg}  (\lambda)  .
$$
Using $\tRic(\eta,\eta) = \tr{\cal{R}}_{\eta}$ and the trace of equation
(\ref{lieetal}) we get 
\begin{equation} \label{trlambda}
\rho' \tr_{\tg} (\lambda) = \rho (\tr_{\tg} ({\cal{L}}_{\eta} \lambda ) 
+ \tr_{\tg} ( \lambda \times \lambda ) ).
\end{equation}
We now start an induction by assuming that 
\begin{equation} \label{assumption}
h =h^{(0)} + \sum_{k=k_0}^{\infty} \frac{x^k}{k!}h^{(k)}
= h^{(0)} + h^{\text{rem}},
\end{equation}
we have already seen that this assumption holds with $k_0=2$.
The Ricci tensor depends smoothly on the metric. Therefore
$\ric - (n-2)h_0 $ vanishes in $x$ to the same order as 
$h^{\text{rem}}$,
that is $ric - (n-2) h = O(x^{k_0})$. If (\ref{assumption}) holds then
$\lambda = O(x^{k_0 -1})$ and 
we conclude using $\rho = x + O(x^3)$ that the lowest order terms in 
(\ref{main}) are the 
terms of order $x^{k_0 -1}$, they are the Lie-derivative term on the left side 
and both terms on the right side. 
Keeping only the lowest order terms in (\ref{main}) and simplifying 
gives us the equation
\begin{equation} \label{fin1}
(n - 1 - k_0)h^{(k_0)} + \tr_{h_0}(h^{(k_0)}) h^{(0)} =0.
\end{equation}
Next we look at the trace $\tr_{h_0}(h^{(k_0)})$.
First we assume that $k_0 > 2$. 
Then if we look at the lowest order terms the 
$\lambda \times \lambda$ term vanishes and after simplification 
\begin{equation} \label{fin2}
(k_0 -2) \tr_{h_0} h^{(k_0)} =0,
\end{equation}
which gives $\tr_{h_0} h^{(k_0)} =0$. If $k_0=2$ we take the 
$h_0$-trace of equation (\ref{fin1}) and get $2(n-2) \tr_{h_0} h^{(2)} =0$,
so $\tr_{h_0} h^{(k_0)}$ vanishes and (\ref{fin1}) becomes
\begin{equation} \label{fin3}
(n - 1 - k_0)h^{(k_0)} =0.
\end{equation}
This means that as long as $k_0 < n-1$, (\ref{assumption}) implies
$h^{(k_0)}=0$, and we have showed that
\begin{equation} \label{seriesresult}
h =h^{(0)} + \sum_{k=n-1}^{\infty} \frac{x^k}{k!}h^{(k)}.
\end{equation}

In order to show that the trace of $A - Id$ vanishes to higher order, we do an
induction on $\tr_{h_0} h^{(k)}$. Assume $\tr_{h_0}h^{(k)} = 0$ for 
$1 \leq k \leq
j$. By the above this holds for $j = n-2$ which is the base case for the
induction. From (\ref{seriesresult}) and (\ref{second}) we have for 
$j <  2n-3$
\begin{eqnarray*}
\tr_{\tg} \lambda &=& \tr_{h_0} \lambda + O(x^{2n-3}) = 
- \half \frac{x^j}{j!} \tr_{h_0} h^{(j+1)} + O(x^{j+1}) ,  \\
\tr_{\tg} ({\cal L}_{\eta} \lambda ) &=& \tr_{h_0} ({\cal L}_{\eta} \lambda ) +
O(x^{2n-4}) = -\half \frac{x^{j-1}}{(j-1)!} \tr_{h_0} h^{(j+1)} + O(x^j) .
\end{eqnarray*}
Now use $n> 3$, 
$\tr_{\tg} (\lambda \times \lambda) = O(x^{2n-4})$ and (\ref{trlambda})
to get, after simplifying and discarding higher order terms, 
$$
 \tr_{h_0} h^{(j+1)} = 0,   \qquad \text{ for } j < 2n-3 .
$$
From the definitions, 
$$
\tg'(A^{-2} X, Y) = \tg(X,Y) = \tg'((Id + H) X,Y),
$$
where $H$ is defined by $\tg'(HX,Y) = h^{\text{rem}}(X,Y)$. Thus from the
above $| H | = O(x^{n-1})$ and $\tr H = O(x^{2n-2})$. 
By definition $A$ and $H$ are self adjoint linear maps w.r.t. the inner
product $\tg'$.
From the equation $A^2 = (Id + H)^{-1}$ we find 
$$
A = \sum_{k=0}^{\infty} \binom{-1/2}{k} H^k = Id - \half H + O(x^{2n-2})
$$
This completes the proof of Lemma \ref{lemma:A}.
~\qed
\begin{remark} Under the assumptions of Fefferman and Graham 
\cite{fefferman:graham}, in case $n$ is even, $h^{(n-1)} =0$, and therefore
instead of (\ref{seriesresult}) we get that $h$ agrees with $h_0$ to any
order in $x$.
~\qed
\end{remark}

\subsection{Proof of Theorem \protect{\ref{thm:conf}}}
\label{sec:proof}
In case $n=3$ $(M,g)$ Einstein implies constant curvature and the theorem
follows from \cite[Theorem 6.9]{andersson:howard:rigid}, so we assume
$n \geq 4$ in the following. \mnote{\lars find some better ref for this}

From equation (\ref{nphi}) we have 
$$
\hn_{X} \phi_0 = (\n_{X} - \bn_{X})\phi_0 - \ihalf(AX-X) \phi_0 ,
$$
outside the compact set where $f \ne 1$.
Lemma \ref{lemma:A} together with (\ref{diffnabla}), (\ref{normgrowth}) 
and the fact that 
$x e^r$ is bounded on $U$ gives
$$
|\hn \phi_0|^2 \leq C (|A^{-1}|^2|\n' A|^2 +|A-Id|^2)x^{-1} = O(x^{2n-3}),
$$
which using  $\mu = \rho^{-n} \widetilde{\mu}$ 
shows that $\hn\phi_0 \in L^2(T^{*}\otimes\Sr)(M)$. 
Clearly $\phi_0 \notin H^1\Sr(M)$ and 
therefore $\phi_0$ is an asymptotic Killing spinor. 

Next we are going to show that the boundary 
integral (where the normal $\nu = \rho\eta$)
\begin{equation}\label{eq:boundterm}
\lim_{x \to 0} \int_{\dM_x} g ((\hn_{\nu}+\nu\hD) \phi_0,\phi_0) d\mu_x 
\end{equation}
vanishes, where $d\mu_x$ is the volume element of the induced metric on
$\dM_x$. 
Introduce an ON frame $\{e'_i\}$ w.r.t. $g'$ adapted to the 
foliation $\dM_x$ so that 
$e'_1 = \nu$ and $e'_i \perp \dM_x$ for $i \geq 2$. Set $e_i = A e'_i$. Then
$\{ e_i\}$ is an ON frame for $g$ adapted to the foliation $\dM_x$. 
In the integrand we have 
$$
(\hn_{\nu}+\nu\hD) \phi_0 = \sum_i (\delta_{1i} + e_1e_i)\hn_{e_i}\phi_0 = 
\sum_i \half \sigma_{1i}\hn_{e_i}\phi_0
$$
where $\sigma_{ij} = [e_i,e_j]$ is the commutator in the Clifford algebra.
We compute
$$
(\n_{X} - \bn_{X})\phi_0 = \quart \sum_{i,j} (\omega_{ij}(X) - 
\overline{\omega}_{ij}(X)) e_i \cdot e_j \cdot \phi_0 =  
\quart \sum_{i,j} (\omega_{ij}(X) - 
\overline{\omega}_{ij}(X)) \sigma_{ij} \cdot \phi_0,
$$
and
\begin{eqnarray} 
g((\hn_{\nu}+\nu\hD) \phi_0,\phi_0) &=&
\frac{1}{8} \sum_{i,j,k} (\omega_{jk}(e_i) - 
\overline{\omega}_{jk}(e_i))g(\sigma_{1i}\sigma_{jk} \phi_0,\phi_0) \nonumber
\\
& & \label{eqn1} -\frac{\bfi}{4} \sum_{i} g(\sigma_{1i}(Ae_i-e_i)\phi_0,\phi_0)
.
\end{eqnarray}
Since the integral in (\ref{eq:boundterm}) is real,
we need not bother about the imaginary terms in (\ref{eqn1}). 
What is left of the first sum is then 
\begin{equation} \label{eqn2}
\sum_{i} (\omega_{i1}(e_i) - \overline{\omega}_{i1}(e_i))g(\phi_0,\phi_0)
+\half \sum_{i,j,k} (\omega_{jk}(e_i) 
-\overline{\omega}_{jk}(e_i))g(\sigma_{1ijk} \phi_0,\phi_0) ,
\end{equation}
where $\sigma_{1ijk} = e_1e_ie_je_k$ if $1,i,j,k$ are all different or
else $\sigma_{1ijk} = 0$.  From (\ref{torsion}) and 
(\ref{difftorsion}) we get the identities
\begin{eqnarray*}
\omega_{i1}(e_i) - \overline{\omega}_{i1}(e_i) &=& 
g(\bT(e_i,e_1),e_i) \\
&=&  g((\n'_{e_1}A)A^{-1}e_i,e_i) - g((\n'_{e_i}A)A^{-1}e_1,e_i)
\end{eqnarray*}
and 
$$
2( \omega_{jk}(e_i) -\overline{\omega}_{jk}(e_i))  =
-g(\bT(e_i,e_j),e_k) + g(\bT(e_i,e_k),e_j) + g(\bT(e_j,e_k),e_i) ,
$$
where the last two terms taken together are symmetric in $ij$ and vanish when 
summed against $\sigma_{1ijk}$. This leaves
\begin{eqnarray} \label{part1}
& &\sum_{i}(g((\n'_{e_1}A)A^{-1}e_i,e_i) - g((\n'_{e_i}A)A^{-1}e_1,e_i))
g(\phi_0,\phi_0) \\
&+& \nonumber \half \sum_{i,j,k} g((\n'_{e_i}A)A^{-1}e_j,e_k)g(\sigma_{1ijk} \phi_0,\phi_0)
\end{eqnarray}
Next we look at the second sum in (\ref{eqn1}). We have 
\begin{eqnarray*}
\sum_{i} g(\sigma_{1i}(Ae_i-e_i)\phi_0,\phi_0) &=&
2\sum_{i} g((\delta_{1i}+e_1e_i)(Ae_i-e_i)\phi_0,\phi_0)\\
&=&2\sum_{i} g(e_1e_i(Ae_i-e_i)\phi_0,\phi_0)\\
&=&2\sum_{i,j} (A_i^j-\delta_i^j)g(e_1e_ie_j\phi_0,\phi_0)
\end{eqnarray*}
where  $A_i^j$ is defined by $A e'_i = A_i^j e'_j$ so that $A e_i = A_i^j e_j$.
\mnote{\lars This gives a contribution in the nonsymmetric case, add this to
remark}
Since $A_i^j$ is symmetric this simplifies to 
$$
-2\sum_{i} (A_i^i-\delta_i^i )g(e_1\phi_0,\phi_0)
=-2(\tr_{g}(A)-n)g(e_1\phi_0,\phi_0)
$$
and (the real part of) the integrand in (\ref{eq:boundterm}) is 
\begin{eqnarray} \label{boundintegrand}
& &\sum_{i}(g((\n'_{e_1}A)A^{-1}e_i,e_i) - g((\n'_{e_i}A)A^{-1}e_1,e_i))
g(\phi_0,\phi_0) \\
&+& \nonumber \half \sum_{i,j,k} g((\n'_{e_i}A)A^{-1}e_j,e_k)g(\sigma_{1ijk} \phi_0,\phi_0)
+ \ihalf (\tr_{g}(A)-n)g(e_1\phi_0,\phi_0) .
\end{eqnarray}

We consider each term separately. To simplify the notation we will
write $f_1 \simeq f_2 $ if $f_1$ and $f_2$ are both $O(x^k)$ for some $k$ and 
$f_1 - f_2 = O(x^{k+n-1})$.
Let $B = A - Id$, note that $B$ is symmetric and $A^{-1} - Id \simeq -B$. In
the following computations we use the summation convention. We have
\begin{eqnarray*}
g(( \n'_{e_1} A ) A^{-1} e_i , e_i ) &=& 
g'( A^{-1} ( \n'_{e_1} A) e'_i, e'_i ) \\
&=& g'( A^{-1} ( \n'_{e_1} (A e'_i) - A \n'_{e_1} e'_i ) , e'_i ) \\
&=& ( e_1 A_i^l )  g'(A^{-1} e'_l , e'_i ) 
+ A_i^l g'(A^{-1} \n'_{e_1} e'_l , e'_i ) 
- g'( \n'_{e_1} e'_i, e'_i ) 
\\
&\simeq& ( e_1 A_i^l )  g'( e'_l , e'_i ) 
+ B_i^l g'( \n'_{e_1} e'_l, e'_i) 
- g'(B \n'_{e_1} e'_i , e'_i )  \\
&=& ( e_1 A_i^i )  
+ B_i^l g'( \n'_{e_1} e'_l, e'_i)
- B_i^l g'(\n'_{e_1} e'_i , e'_l ) \\
&=& ( e_1 A_i^i )  ,
\end{eqnarray*}
where in the last row we used the fact that $B$ is symmetric. Similarly we
have
\mnote{the last term with $B$ is there in the general expression for the mass}
\begin{eqnarray*}
g(( \n'_{e_i} A ) A^{-1} e_1 , e_i ) &=&
g'( A^{-1} ( \n'_{e_i} A) e'_1, e'_i ) \\
&=& g'( A^{-1}  (Id - A) \n'_{e_i}  e'_1 , e'_i ) \\
&\simeq& - g'(B \n'_{e_i}  e'_1 , e'_i ) \\
&=& - B_i^l g'( \n'_{e_i} e'_1 , e'_l ) \\
&\simeq& - B_i^l g'( \n'_{e'_i} e'_1, e'_l)  \\
&=& - B_i^l ( \rho'  \delta_{i1} \delta_{1l} 
+ \rho \tg'( \n'_{\te'_i} \te'_1, \te'_l ) ) , 
\end{eqnarray*}
where $\te'_i = \rho^{-1} e'_i$ form an ON frame for $\tg' = \rho^2 g'$.
We now use the formula 
$$
\tg'(\n'_{X} \te'_i, \te'_j) = 
\omega'_{ij}(X) = \tomega'_{ij}(X) - \rho^{-1} ( (\te'_j \rho) \tg'( \te'_i,X)
- (\te'_i \rho) \tg'( \te'_j,X) )
$$
for the connection coefficients for $\n'$ w.r.t. the frame $\{ \te'_i \}$.
Note that $\tomega'_{1l}(\te'_i) = 0$ for the cylindrical metric. 
Therefore 
$$
\rho \tg'( \n'_{\te'_i} \te'_1, \te'_l ) = - \rho' \delta_{1l} \delta_{1i} 
+ \rho' \delta_{il}
$$
and hence, 
$$
g(( \n'_{e_i} A ) A^{-1} e_1 , e_i ) \simeq - \rho' B_i^i .
$$
\newcommand{\sis}{\sigma_{1ijk}}
Next we consider the terms in (\ref{boundintegrand}) containing $\sis$. We
have
\begin{eqnarray*}
g(( \n'_{e_i} A ) A^{-1} e_j , e_k ) \sis &=&
 g'( A^{-1} ( \n'_{e_i} A) e'_j, e'_k )\sis  \\ 
&=&  g'( A^{-1} ( \n'_{e_i} (A e'_j) - A \n'_{e_i} e'_j) , e'_k )  \sis \\
&=& \left ( ( e_i A_j^l )  g'(A^{-1} e'_l , e'_k )
+ A_j^l g'(A^{-1} \n'_{e_i} e'_l , e'_k ) \right. \\
&& \hskip 0.2in \left.
- g'( \n'_{e_i} e'_j, e'_k ) \right ) \sis \\
&\simeq&
\left ( 
e_i A_j^k  + B_j^l g'(\n'_{e_i} e'_l, e'_k ) 
- B_k^l g'( \n'_{e_i} e'_j , e'_l) \right ) \sis \\
&=& 
\left ( 
e_i A_j^k  + B_j^l g'(\n'_{e_i} e'_l, e'_k ) 
+ B_k^l g'( \n'_{e_i} e'_l , e'_j) \right ) \sis \\
&=& 0 . 
\end{eqnarray*}
This means that $Q = g(( \n'_{e_i} A ) A^{-1} e_j , e_k ) \sis$ satisfies  
$Q = O(x^{2n-2})$.
Finally, note that $\tr_g A - n \simeq B_i^i $ where the components are
defined w.r.t. the frame $\{ e'_i \}$.
This gives the following expression for the mass 
\begin{eqnarray}
\lim_{x \to 0} \int_{\dM_x} g ((\hn_{\nu}+\nu\hD) \phi_0,\phi_0) d\mu_x 
&=& \lim_{x \to 0} \int_{\dM_x}  \nonumber 
\bigg ( \big ( (  e_1 A_i^i ) + \rho' B_i^i  \big ) g(\phi_0.\phi_0)  \\
&&  \hskip -0.5in + \half g(Q \phi_0, \phi_0) 
+ \ihalf B_i^i g(e_1 \phi_0, \phi_0) \bigg ) 
d\mu_x  . \label{massexpr}
\end{eqnarray}
From Lemma \ref{lemma:A} we have $e_1 A_i^i = O(x^{2n-2}), B^i_i =
O(x^{2n-2})$ and as we have seen, $Q = O(x^{2n-2})$. 
Further, using (\ref{normgrowth}) and $x e^r$ bounded on $U$
together with 
$d\mu_x = \rho^{-(n-1)} d\mu(h)$ we find that the integrand in (\ref{massexpr})
is $O(x^{n-2})$ which vanishes as $x \to 0$ and hence the mass is zero.

It follows from Lemma \ref{lemma:basic}
that there is a Killing spinor $\phi$ for $g$ associated to any Killing spinor
$\phi_0$ for the hyperbolic background metric $g'$ and hence there is a full
set of linearly independent Killing spinors for $g$. This proves that 
\mnote{\lars discuss linear independence!}
$\hR =0$ and hence $(M,g)$ has sectional curvature $-1$ and is isometric to
hyperbolic space $\hy$.
~\qed

\mnote{\lars here is a remark on the form of the mass for $A$ general
commented out}

\subsection{Spherical space form boundary}
\begin{thm}
\label{thm:conf2}
Let $M$ conformally compact Einstein spin manifold with conformal boun\-dary $\dM = S^{n-1}/\Gamma$.
Suppose that the spin structure on a collar neighborhood $U$ of $\dM$ is 
equivalent to the spin structure on $\hp /\Gamma$ defined by the lift $\tG$ of $\Gamma$.
If $\tG$ fixes some non--zero spinor $u \in \Sr$ then $ \Gamma = \{ 1 \}$ and $M$ is isometric to $\hy$.
\end{thm}
\proof
In case $n=3$ $(M,g)$ Einstein implies constant curvature and the theorem
follows from \cite[Theorem 6.9]{andersson:howard:rigid}, so we assume
$n \geq 4$ in the following. \mnote{\lars find some better ref for this}

By passing to a cover of $U$, the estimates of Lemma \ref{lemma:A} 
apply without change. 
From the assumptions it follows that $U$ with the hyperbolic background metric $g'$ and its induced
spin structure has a Killing spinor. Therefore by the same argument as in the proof of Theorem 
\ref{thm:conf} there is an asymptotic Killing spinor with mass zero. 

By Lemma \ref{lemma:basic} there is a Killing spinor on $(M,g)$ and therefore $(M,g)$ is a warped
product as in Theorem \ref{thm:baum}. 

It follows from the estimates of Lemma \ref{lemma:A} that an integral of the form 
$$
\int_{\dM_x} (|{\n'}^2A| + |\n' A|^2 + |A-Id|) d\mu_{x}
$$
is bounded. 
\mnote{discuss $r,t,x,....$}
Therefore since we are considering only $n \geq 4$ the proof of Lemma 
\ref{lemma:done} applies to the present case.  
~\qed

\section{Four dimensions}
\label{sec:4d}
We now consider the four dimensional case. Here we have isomorphisms
$$ 
\Si(4) = \SU(2) \times \SU(2) = \Symp(1) \times \Symp(1) .
$$ 
\mnote{\lars change here, added $\Co^{2*}$}
The projection $\Si(4) \to \SO(4)$ takes $(p,q)$ to the map $v @>>> pvq^*$ 
(quaternion multiplication). The spinor representation is 
$\SU(2) \times \SU(2)$ acting on $\Sr= \Sr^+ \oplus \Sr^- = \Co^2 \oplus \Co^{2*} $. 

\subsection{ALH ends in four dimensions} Suppose that 
$\Gamma \subset \SO(4)$ has a lift to $\tG \subset \Si(4)$ 
such that for all $\tilde{\gamma} \in \tG$, $\tilde{\gamma} \cdot u = u$
where $0 \neq u \in \Sr$. 
Then the same holds for the parts $u_+,u_-$, by a choice of orientation 
we may assume $u_- \neq 0$. This means that for all $\tilde{\gamma}$ 
the part in the second $\SU(2)$ factor, $\tilde{\gamma}_-$, has one 
eigenvalue equal to one and since the determinant is one we must have 
\mnote{\lars why determinant one?}
$\tilde{\gamma}_- = Id$. So $\Sr^-$ is fixed by $\tG$ and $\tG$ is a 
subgroup of the first $\SU(2)$ factor, which also means 
that $\Gamma \subset \Symp(1) = \SU(2) \subset \SO(4)$. The same 
reasoning gives that every finite subgroup of $\SU(2)$ acts freely on the 
sphere. We conclude;
\begin{prop}
\mnote{\lars What about higher dimensional case?}
If $\Gamma$ is a finite subgroup of $\SU(2)$ then with the above choice of
orientation, the spinors $\phi_u$ with 
$u \in \Sr^-$ give Killing spinors on the quotient $({\Bbb H}_*^4) / \Gamma$ 
provided we choose 
the spin structure 
$$
\Si(({\Bbb H}_*^4) / \Gamma) = \Si({\Bbb H}_*^4) / \tG
$$ 
defined by the lift $\gamma @>>> \tilde{\gamma} = (\gamma,Id)$. 
Except for reversing orientation these are the only cases allowing 
Killing spinors.
\end{prop}
So we can only find asymptotic Killing spinors on an asymptotically 
locally hyperbolic four manifold if the fundamental group of the 
locally hyperbolic end
is a finite subgroup of $\SU(2)$. The following groups are up to 
conjugation all finite subgroups of $\SU(2)$ (\cite[Thm. 2.6.7]{wolf}). 
\begin{itemize}\closeup
\item[${\mathbf A}_{\mathbf n}$:] The cyclic group of order $n$ generated by
$
z = \begin{pmatrix} \zeta & 0 \\ 0 & \zeta^{-1} \end{pmatrix} 
$
where $\zeta = e^{2\pi i/n}$.
\item[${\mathbf D}_{\mathbf n}^{\mathbf *}$:] The binary dihedral group of order 4n, this consists
of $\{z^a,jz^a\}^{2n-1}_{a=0}$ where 
$
z = \begin{pmatrix} \zeta & 0 \\ 0 & \zeta^{-1} \end{pmatrix} 
$
with $\zeta = e^{2\pi i/2n}$ and
$
j = \begin{pmatrix} 0 & 1 \\ -1 & 0 \end{pmatrix}. 
$
\item[${\mathbf T}^{\mathbf *}$:] The binary tetrahedral group.
\item[${\mathbf O}^{\mathbf *}$:] The binary octahedral group.
\item[${\mathbf I}^{\mathbf *}$:] The binary icosahedral group.
\end{itemize} 
For these groups there is a canonical lift to $\Si(4)$ and an associated
canonical spin-structure on the end which we call the trivial spin structure.
\mnote{\lars is it correct to call this trivial?}
Any other lift of $\Gamma$ to $\Si(4)$ is given up to conjugation by 
an element $\kappa \in \operatorname{Hom}(\Gamma,\Za_2)$ as follows
$$
\Gamma \ni \gamma @>>> \tilde{\gamma} = \kappa(\gamma)(\gamma,Id) \in\Si(4).
$$ 
The elements of $\operatorname{Hom}(\Gamma,\Za_2)$ are
\begin{itemize}\closeup
\item[${\mathbf A}_{\mathbf n}$:] $\kappa_0 =1$ and if $n$ is even $\kappa_1$ defined 
by $\kappa_1(z) = -1$.
\item[${\mathbf D}_{\mathbf n}^{\mathbf *}$:] $\kappa_{pq}$, $p,q =0,1$ defined by 
$\kappa_{pq}(z) = (-1)^p$ and $\kappa_{pq}(j) = (-1)^q$.
\item[${\mathbf T}^{\mathbf *}$:] 1.
\item[${\mathbf O}^{\mathbf *}$:] $\kappa_0 = 1$ and  $\kappa_1$ which is nontrivial.
\item[${\mathbf I}^{\mathbf *}$:] 1.
\end{itemize} 
These are thus also the possible spin structures on an asymptotically locally
hyperbolic end. 
\mnote{\lars added locally}
\subsection{Detecting the spin structure on the end}
We will now see that the spin structure on an asymptotically locally 
hyperbolic end can almost be detected by the signature of the manifold. For 
simplicity we assume for the rest of this section that $M$ is an 
asymptotically locally hyperbolic manifold with only one end.
The problem of detecting the spin-structure is a purely topological one so
we can forget the original metric on $M$ and instead put on $M$ a metric
which is a product $dt^2 + g_1$ on the end $(0,1]\times (S^3/\Gamma)$,
where $g_1$ is the spherical metric. Note that we have added a boundary
isometric to the spherical space form $(S^3/\Gamma)$. We can then apply the 
Atiyah--Patodi--Singer index theorem to compute the index of the Dirac 
operator and the signature of the manifold
\mnote{\lars add remark on harmonic spinors on $\dM$}
\begin{equation}
\operatorname{ind}(D) = -\frac{1}{24}\int_M p_1 - \eta_D(S^3/\Gamma),
\end{equation} 
\begin{equation}
\sigma(M) = \frac{1}{3}\int_M p_1 - \eta_{\sigma}(S^3/\Gamma).
\end{equation}
where $\eta_{\sigma}(S^3/\Gamma)$ and $\eta_D(S^3/\Gamma)$ are 
the eta-invariants of the Signature- and the Dirac-operator on the 
boundary. These formulas combine to 
$$ 
\sigma(M) + 8\operatorname{ind}(D)=-\eta_{\sigma}(S^3/\Gamma)-8\eta_D(S^3/\Gamma).
$$
On a four dimensional manifold the index of the Dirac operator is 
even since the spin representation is quaternionic, this means that 
the $8\operatorname{ind}(D)$ term vanishes modulo 16 and
\begin{equation} \label{sigmarel}
\sigma(M) \equiv -\eta_{\sigma}(S^3/\Gamma)-8\eta_D(S^3/\Gamma) 
\quad (\mod 16). 
\end{equation} 
The right hand side of (\ref{sigmarel}) is known as the Rochlin invariant. 
This expression is useful since the eta invariants are 
explicitly computable \cite{gilkey:spherical} and $\eta_D$ involves the 
spin structure via $\kappa$.
The details of the computation of the eta invariants will be described in
\cite{andersson:dahl:future}. We summarize the results for the case of 
$\Gamma \subset \SU(2)$ in Table \ref{tab:tab1}.
From (\ref{sigmarel}) we now have 
\begin{thm} \label{thm:signature} Let $M$ be a spin manifold of dimension $4$ 
with boundary $\dM = S^3/\Gamma$ where 
$\Gamma \subset \SU(2)$. Fix the orientation of $M$ so that the induced
orientation on $\partial M$ is standard.
Then if $\sigma(M)(\mod 16)$  does not take the value 
corresponding to  one of the nontrivial spin structures in Table \ref{tab:tab1},
then the induced spin structure on $\dM$ is trivial.
\end{thm}
In the following we assume the orientation of $M$ is chosen as in Theorem
\ref{thm:signature}.
Combining Theorem \ref{thm:signature} with Theorem \ref{thm:main1} proves
\begin{cor}\label{cor:4dALH} Let $(M,g)$ be a $4$ dimensional complete spin manifold with one ALH
end 
with group $\Gamma \subset \SU(2)$ and $\hs \geq 0$. 
Then $\sigma(M)(\mod 16)$  does not take the
value of the trivial spin structure corresponding to $\Gamma$ in 
Table  \ref{tab:tab1} unless $(M,g)$ is isometric to ${\Bbb H}^4$. Therefore
the allowed groups and signatures are those listed in Table \ref{tab:tab2}.
\end{cor}
 
Similarly, Theorem \ref{thm:signature} combined with Theorem \ref{thm:conf2}
proves
\begin{cor}\label{cor:4dEin} Let $(M,g)$ be a $4$ dimensional Einstein spin 
manifold which 
is conformally compact with conformal boundary $\dM = S^3/\Gamma$ with $\Gamma \subset \SU(2)$.
Then $\sigma(M)(\mod 16)$  does not take the
value of the trivial spin structure corresponding to $\Gamma$ in Table
\ref{tab:tab1} unless $(M,g)$ is isometric to ${\Bbb H}^4$.
Therefore 
the allowed groups and signatures are those listed in Table \ref{tab:tab2}.
\end{cor}
\mnote{discuss counterexamples a'la LeBrun}
\begin{table}[htp]
\caption[]{Allowed groups and signature}
\label{tab:tab2}
\begin{center}
\begin{tabular}{l|c}
Group			&  $\sigma(M) (\mod 16)$  \\
\hline 
$\{ id \}$ 		& 0    \\
\hline
${\mathbf A_n}$ ($n$ even) 	& 1    \\
\hline
${\mathbf D_n^*}$  		& -n,-2,0      \\
\hline
${\mathbf O^*} $     	&    $-1$   \\
\end{tabular}
\end{center}
\end{table}

\def\arraystretch{1.6}
\begin{table}[htp]
\caption[]{Eta invariants and signature}
\label{tab:tab1}
\begin{center}
\begin{tabular}{l| l| c| c| c}
Group			& Spin structure & $ \eta_{\sigma}$ & $\eta_D$ & $\sigma(M) (\mod 16)$  \\
\hline
${\mathbf A_n}$ ($n$ even) 	& $\kappa_0$     & $\frac{(n-1)(n-2)}{3n}$ &$\frac{n^2-1}{12n}$&  $-n+1$    \\
                       	& $\kappa_1$     &                       &$-\frac{n^2+2}{12n}$&    $1$      \\
\hline 
${\mathbf A_n}$ ($n$ odd) 	& $\kappa_0$     & $\frac{(n-1)(n-2)}{3n}$ &$\frac{n^2-1}{12n}$&   $-n+1$   \\
\hline
${\mathbf D_n^*}$  		& $\kappa_{00}$  & $\frac{2n^2+1}{6n}$ &$\frac{4n^2+12n-1}{48n}$ &  $-n-2$      \\
               		& $\kappa_{01}$  &                   &$\frac{4n^2-1}{48n}$ &   $-n$         \\
               		& $\kappa_{10}$  &                 &$-\frac{2n^2-12n+1}{48n}$ &   $-2$     \\

               		& $\kappa_{11}$  &                 &$-\frac{2n^2+1}{48n}$ &   $0$           \\
\hline
${\mathbf T^*}$      	&                & $\frac{49}{36}$ &$\frac{167}{288}$ &   $-6$     \\
\hline
${\mathbf O^*} $     	& $\kappa_0$     & $\frac{121}{72}$ &$\frac{383}{576}$ &    $-7$   \\
               		& $\kappa_1$     &                  &$-\frac{49}{576}$ &    $-1$    \\
\hline
${\mathbf I^*} $     	&                & $\frac{361}{180}$ &$\frac{1079}{1440}$ &   $-8$     \\
\end{tabular}
\end{center}
\end{table}
\def\arraystretch{1}


\ifx\undefined\bysame
\newcommand{\bysame}{\leavevmode\hbox to3em{\hrulefill}\,}
\fi

\end{document}